\documentclass[letterpaper,10pt,conference]{ieeeconf} 

\IEEEoverridecommandlockouts 

\usepackage{amsmath,amssymb}
\usepackage{mathtools}
\usepackage{array}
\usepackage{booktabs}
\usepackage{hyperref}
\usepackage[compress]{cite}
\usepackage{subcaption}
\usepackage{tikz,pgfplots}
\usetikzlibrary{positioning}
\usepackage[linesnumbered,ruled,vlined]{algorithm2e}

\title{\LARGE\bf%
Online Complexity Estimation for Repetitive Scenario Design%
}

\author{Guillaume O.~Berger$^{1}$, and Rapha\"el M.~Jungers$^{1}$
\thanks{$^{1}$Both authors are with ICTEAM, UCLouvain, Louvain-la-Neuve, Belgium.
GB is an FNRS Postdoctoral Researcher.
RJ is an FNRS honorary Research Associate.
This project has received funding from the European Research Council (ERC) under the European Union's Horizon 2020 research and innovation programme under grant agreement no.~864017 (L2C), from the Horizon Europe programme under grant agreement no.~101177842 (Unimaas), and from the ARC (French Community of Belgium) (SIDDARTA).
Emails: {\{guillaume.berger,raphael.jungers\}@uclouvain.be}.}%
}

\newcommand{\R}{\mathbb{R}}
\newcommand{\N}{\mathbb{N}}
\DeclareMathOperator*{\argmax}{arg\,max}

\newcommand{\Xsf}{\mathsf{X}}
\newcommand{\Prob}{\mathsf{P}}
\newcommand{\Exp}{\mathsf{E}}

\newcommand{\dd}{\mathrm{d}}
\newcommand{\calG}{\mathcal{G}}
\newcommand{\calD}{\mathcal{D}}
\newcommand{\Alg}{\mathsf{A}}

\newcommand{\xt}{\tilde{x}}
\newcommand{\calN}{\mathcal{N}}

\newcommand{\Bsf}{\mathsf{B}}
\newcommand{\gb}{\mathbf{g}}
\newcommand{\vh}{\hat{v}}
\newcommand{\Nh}{\bar{n}}
\newcommand{\oneb}{\boldsymbol{1}}
\newcommand{\Nmax}{N_{\mathrm{max}}}
\newcommand{\thetat}{\tilde\theta}
\newcommand{\gbt}{\tilde\gb}

\newcommand{\newgb}{\color{blue}}
\renewcommand{\newgb}{}

\newtheorem{definition}{Definition}
\newtheorem{theorem}{Theorem}
\newtheorem{lemma}{Lemma}
\newtheorem{proposition}{Proposition}
\newtheorem{corollary}{Corollary}
\newtheorem{problem}{Problem}
\newtheorem{remark}{Remark}
\newtheorem{example}{Example}

\newif\ifextended
\extendedtrue

\begin{document}

\maketitle
\thispagestyle{empty}
\pagestyle{empty}

\begin{abstract}
We consider the problem of repetitive scenario design where one has to solve repeatedly a scenario design problem and can adjust the sample size (number of scenarios) to obtain a desired level of risk (constraint violation probability).
We propose an approach to learn on the fly the optimal sample size based on observed data consisting in previous scenario solutions and their risk level.
Our approach consists in learning a function that represents the pdf (probability density function) of the risk as a function of the sample size.
Once this function is known, retrieving the optimal sample size is straightforward.
We prove the soundness and convergence of our approach to obtain the optimal sample size for the class of fixed-complexity scenario problems, which generalizes fully-supported convex scenario programs that have been studied extensively in the scenario optimization literature.
We also demonstrate the practical efficiency of our approach on a series of challenging repetitive scenario design problems, including non-fixed-complexity problems, nonconvex constraints and time-varying distributions.
\end{abstract}


\section{INTRODUCTION}\label{sec:intro}

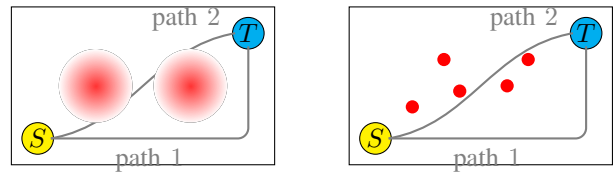
\begin{figure}[t]
    \centering
    \begin{subfigure}[b]{0.48\linewidth}
        \centering
        \begin{tikzpicture}[scale=0.7]
            \draw[black] (0,0) rectangle (5,3);
            \draw[black,fill=yellow] (0.5,0.5) circle (0.3) node {$S$};
            \draw[black,fill=cyan] (4.5,2.5) circle (0.3) node {$T$};
            \draw[gray,thick,rounded corners] (0.75,0.5) -- (4.5,0.5) node[pos=0.5,below] {path $1$} -- (4.5,2.25);
            \draw[gray,thick,rounded corners] (0.75,0.5) to[out=10,in=190] node[pos=0.8,above left=1mm and -5mm] {path $2$} (4.25,2.5);
            \fill[inner color=red!80!white,outer color=white,draw=none] (1.6,1.5) circle (0.7);
            \fill[inner color=red!80!white,outer color=white,draw=none] (3.4,1.5) circle (0.7);
        \end{tikzpicture}
        \caption{Model-based problem: pdf of the position of the obstacles.}
        \label{fig:example-path-distribution}
    \end{subfigure}
    \hfill
    \begin{subfigure}[b]{0.48\linewidth}
        \centering
        \begin{tikzpicture}[scale=0.7]
            \draw[black] (0,0) rectangle (5,3);
            \draw[black,fill=yellow] (0.5,0.5) circle (0.3) node {$S$};
            \draw[black,fill=cyan] (4.5,2.5) circle (0.3) node {$T$};
            \draw[gray,thick,rounded corners] (0.75,0.5) -- (4.5,0.5) node[pos=0.5,below] {path $1$} -- (4.5,2.25);
            \draw[gray,thick,rounded corners] (0.75,0.5) to[out=10,in=190] node[pos=0.8,above left=1mm and -5mm] {path $2$} (4.25,2.5);
            \node[fill=red,circle,minimum width=5pt,inner sep=0pt] (x) at (1.8,2) {};
            \node[fill=red,circle,minimum width=5pt,inner sep=0pt] at (2.1,1.4) {};
            \node[fill=red,circle,minimum width=5pt,inner sep=0pt] at (1.2,1.1) {};
            \node[fill=red,circle,minimum width=5pt,inner sep=0pt] at (3.4,2) {};
            \node[fill=red,circle,minimum width=5pt,inner sep=0pt] at (3,1.5) {};
        \end{tikzpicture}
        \caption{Scenario problem: sampled positions of the obstacles.}
        \label{fig:example-path-sample}
    \end{subfigure}
    \caption{Comparison between the model-based problem (a) and the scenario problem (b).
    The scenario design algorithm will choose path 2 over path 1, because both paths avoid the obstacles in all sampled positions (b) but path 1 is shorter.
    We see (a) that path 1 avoids the obstacles with probability 1, while path 2 avoids the obstacles with probability $1-\epsilon$, with $\epsilon>0$.
    Under some conditions on the problem, and if enough samples are used, this $\epsilon$ can be arbitrarily small, with high confidence.}
    \label{fig:example-path}
\end{figure}

This paper is concerned with the problem of using the optimal amount of data in repetitive scenario design.
Scenario design is a powerful tool for designing ``optimal'' solutions while guaranteeing that some random constraint $g(x)\leq0$, where $g:\Xsf\to\R$ is picked randomly, is satisfied with high probability; see, e.g.,~\cite{calafiore2006thescenario,campi2008theexact,calafiore2010random,campi2011asamplinganddiscarding,margellos2014ontheroad,esfahani2015performance,grammatico2016ascenario,campi2018waitandjudge,campi2018ageneral,yang2019chanceconstrained,rocchetta2021ascenario,romao2021tight,garatti2022risk,romao2023ontheexact,campi2023compression}.
The principle of scenario design is to replace the random constraint by $N$ i.i.d.~samples (called \emph{scenarios}) of it, namely $g_i(x)\leq0$ for $i=1,\ldots,N$.
The final solution is then selected---e.g., by optimizing some preference criterion $J(x)$---among all solutions that satisfy the sampled constraints.
To fix ideas, let us consider the example of optimal path planning in an uncertain environment consisting of randomly positioned obstacles (see Fig.~\ref{fig:example-path}).
In this context, the constraint $g(x)\leq0$ is to avoid the obstacles and its probability distribution is given by the probability distribution of the position of the obstacles (see Fig.~\ref{fig:example-path-distribution}).
The scenario design approach consists in drawing $N$ sampled positions of the obstacles, and finding a path that avoids the obstacles in all---or a predefined fraction---of the sampled positions (see Fig.~\ref{fig:example-path-sample}).

Under some conditions on the problem, and if $N$ is large enough, probabilistic guarantees on the constraint violation probability (called the \emph{risk}) of the solution returned by the scenario design algorithm can be given by the theory of scenario design~\cite{calafiore2006thescenario,campi2018ageneral,alamo2009randomized,calafiore2010random,campi2011asamplinganddiscarding,margellos2015ontheconnection}.
However, a large sample size $N$ can severely affect the cost of finding a feasible solution to the sampled problem.
For instance, if the constraints are nonconvex (as in Fig.~\ref{fig:example-path} for instance), the cost of solving the sampled problem can be exponential in $N$.
Moreover, in cases where one can tolerate some probability of failure, sampling too many constraints can lead to overly conservative solutions~\cite{campi2011asamplinganddiscarding}.
Therefore, finding the smallest $N$ that guarantees an upper bound on the risk is of paramount importance for practical applications.

We address this question in the context of \emph{repetitive} scenario design.
Repetitive scenario design is when {\newgb a similar or slowly-varying} scenario design task is performed repeatedly.
In other words, it is an online version of scenario design that accounts for continued updates in the data and the environment{\newgb, assuming small magnitude of the updates}.
Online path planning is a good example, where path planning (Fig.~\ref{fig:example-path}) is performed repeatedly to account for changes in the position of the obstacles {\newgb and the initial condition of the system}.
Because the task is repeated{\newgb, and under sufficiently slow variations,} we can use the information gathered from previous computations to optimize the sample size $N_t$ in future steps, so that ultimately the sample size is optimal, i.e., $N_t\to N_\star(t)$, where $N_\star(t)$ is the optimal sample size at time $t$.
Nevertheless, this must be done in a cautious way because we cannot {\newgb exceed the risk tolerance} too often during the task.
This precludes naive strategies like ``augmenting $N_t$ if we exceed the risk tolerance at step $t-1$, and decreasing $N_t$ if we meet the risk tolerance at step $t-1$''.

Instead, we propose an approach consisting in learning a function $f_\theta(v,N)$---with arguments $v\in[0,1]$ and $N\in\N$, and parameter $\theta$---that aims to approximate for each $N\in\N$ the probability density function (pdf) of the risk of the solution returned by the scenario design algorithm with $N$ i.i.d.~samples.
To fit $\theta$, we use previously collected data, which consists of triples $(N_s,x_s)$, for $s=1,\ldots,t$, where $t$ is the current step, $x_s$ is the solution computed at step $s$.
The advantage of this approach is that{\newgb, through $\theta$,} we can retrieve information about the risk pdf $f_\theta(\cdot,N)$, even for values of $N$ that have never been used, and in this way infer the optimal sample size $N_\star$ without exceeding the risk tolerance $\epsilon$ too often (Section~\ref{sec:approach}).
We demonstrate the correctness and convergence (to the optimal sample size) of our approach for the class of ``fixed-complexity'' scenario problems, when $f_\theta(\cdot,N)$ is the pdf of the \emph{Beta distribution} with parameters $\theta$ and $N-\theta+1$ (Section~\ref{sec:theory}).
{\newgb The fixed-complexity scenario problems---which generalize fully-supported random convex programs~\cite{campi2008theexact,calafiore2010random}---serve in this paper as a paradigmatic class of slowly-varying scenario problems.}
We demonstrate the practical applicability and efficiency of our approach on a wide range of repetitive scenario design problems, including non-fixed-complexity scenario problems, nonconvex constraints, time-varying distribution, etc.~(Section~\ref{sec:experiments}).

\subsection*{Comparison with the literature}

The seminal works~\cite{campi2008theexact,calafiore2010random} on convex scenario optimization introduced the notion of fully-supported problems.
For these problems, the pdf of the risk is known to be proportional to $v^{d-1}(1-v)^{N-d}$ (Beta distribution), where $d$ is the dimension of the decision variable $x$.
For general convex optimization problems, only an upper bound on the cdf of the risk is known~\cite{littlestone1986relating,campi2008theexact,calafiore2010random,margellos2015ontheconnection}.
This upper bound can be quite conservative, especially when the number of ``support constraints'' is small compared to the dimension of the decision variable~\cite{garatti2022risk}.
One solution to this problem is the ``wait-and-judge'' scenario approach~\cite{campi2018waitandjudge,garatti2021therisk}, where the ``complexity of the set of sampled constraints'' is computed after the samples are drawn in order to derive an upper bound on the cdf of the risk.
Yet, computing the complexity of a set of constraints can be challenging (typically, exponential in $N$) and the resulting bound on the risk is still conservative~\cite{garatti2021therisk}.
Our approach, by contrast, tries to learn directly the pdf of the risk as a function of $N$, thereby removing most of the conservatism.
For the class of fixed-complexity convex programs studied in Section~\ref{sec:approach}, the parameter $\theta$ represents the complexity of the problem {\newgb and is assumed to be fixed (as an ideal case of slowly varying).}
In this case, the approach can be seen as an indirect way of learning this complexity, not requiring computations of the complexity of sample sets.

The recent papers~\cite{calafiore2016repetitive,garatti2023complexity} also use repeated calls to a scenario design algorithm in order to solve a design problem with random constraints.
However, their objective is different: they aim to solve the problem only once (not continually as us).
For that, they increase progressively the number $N$ of samples until a suitable solution is found.
These approaches are not well suited for our online scenario design framework because there is no mechanism for decreasing $N$, so that $N$ will be conservative with probability one in the long run.
Furthermore, the approach in~\cite{garatti2023complexity} requires to compute the complexity of sample sets, which is something that we want to avoid.

One limitation of our approach is that we need to evaluate the risk $\Prob[g(x_t)>0]$ (cf.~Section~\ref{ssec:data-collection}) of the solution $x_t$ at each step $t$.
Computing the risk exactly can be costly in some applications.
Alternatively, one can use sample-based approximation methods as in\ifextended\ Appendix~\ref{app:risk-evaluation}\else~\cite[Appendix~A]{berger2025online}\fi; see also~\cite[Section~II.A]{calafiore2016repetitive}.
Clearly, using sample-based methods for evaluating the risk requires to use extra samples; this additional sample complexity, as well as a way to optimize it {\newgb(e.g., by reusing previous samples or tuning the accuracy of the approximation)}, is not included in our analysis (we leave it for future work).
The results of this paper remain nevertheless valuable also in contexts where risk is evaluated via sampling.
For example, when the sampling cost is negligible compared to the computational cost of solving large-scale scenario design problems (e.g., nonconvex problems with complexity exponential in  $N$), it becomes worthwhile to invest effort in accurate risk evaluation if doing so reduces the size of the scenario design problems to be solved.

\ifextended\relax
\else
All proofs are available in the extended version~\cite{berger2025online}.
\fi

\paragraph*{Notation}

$\N$ is the set of nonnegative integers.
For $n\in\N$, we let $[n]=\{1,\ldots,n\}$.
$\lVert\cdot\rVert_2$ is the Euclidean norm, and $\lVert\cdot\rVert_\infty$ is the $L^\infty$-norm.
The following functions will be useful:
\begin{itemize}
    \item The \emph{Beta function}, defined for all $a,b>0$ by $B(a,b)=\int_0^1 v^{a-1}(1-v)^{b-1}\,\dd v$.
    \item The \emph{Gamma function}, defined for all $z>0$ by $\Gamma(z)=\int_0^\infty t^{z-1}e^{-z}\,\dd t$.
    \item The \emph{Digamma function}, defined for all $z>0$ by $\Psi(z)=\Gamma'(z)/\Gamma(z)$~\cite[\S6]{abramowitz1972handbook}.
\end{itemize}\vskip0pt

\section{PROBLEM STATEMENT}\label{sec:problem}

We consider an optimization problem of the form
\begin{equation}\label{eq:optim-robust}
\min_{x\in\Xsf}\; J(x,\xi) \quad \text{s.t.} \quad g(x,\xi)\leq0 \quad \forall\,g\in\calG,
\end{equation}
where $J:\Xsf\times\Xi\to\R$, $\calG$ is a set of functions from $\Xsf\times\Xi$ to $\R${\newgb, and $\xi\in\Xi$ is an external parameter}.
\eqref{eq:optim-robust} is called the \emph{robust design problem} because the constraint $g(x,\xi)\leq0$ must be satisfied for \emph{all} $g\in\calG$.
A relaxation---and sometimes more realistic version---of the robust problem is the \emph{chance-constrained design problem}:
\begin{equation}\label{eq:optim-cc}
\min_{x\in\Xsf}\; J(x,\xi) \quad \text{s.t.} \quad \Prob[g(x,\xi)>0]\leq\epsilon,
\end{equation}
where $\Prob$ is a probability measure on $\calG$, and $\epsilon\in[0,1]$ is an upper bound on the \emph{risk}, i.e., the probability of violating the constraint $g(x,\xi)\leq0$.
The difference between~\eqref{eq:optim-cc} and~\eqref{eq:optim-robust} is that a feasible solution $x$ of~\eqref{eq:optim-cc} is allowed to violate the constraint $g(x,\xi)\leq0$ for some values of $g\in\calG$, provided the probability measure of these values does not exceed $\epsilon$.
\eqref{eq:optim-cc} can be very challenging to solve, as it typically involves a nonconvex optimization problem~\cite{charnes1959chanceconstrained}.
A way to circumvent this is to consider instead the \emph{scenario design problem}:
\begin{equation}\label{eq:optim-scenario}
\min_{x\in\Xsf}\; J(x,\xi) \quad \text{s.t.} \quad g_i(x,\xi)\leq0 \quad \forall\,i\in[N],
\end{equation}
where for each $i\in[N]$, $g_i\in\calG$.
Clearly,~\eqref{eq:optim-scenario} is a relaxation of~\eqref{eq:optim-robust} since only a subset of $\calG$ is used.
For this reason, a feasible solution of~\eqref{eq:optim-scenario} is not expected to be feasible for~\eqref{eq:optim-robust} in general.
However, if $N$ is large and the $g_i$'s are sampled i.i.d.~according to $\Prob$, one can expect that the solution of~\eqref{eq:optim-scenario} is feasible for~\eqref{eq:optim-cc}.
Note that when the $g_i$'s are sampled at random, the solution of~\eqref{eq:optim-scenario} is also a random variable (since it depends on the samples).
Therefore, the property that the solution of~\eqref{eq:optim-scenario} is feasible for~\eqref{eq:optim-cc} is a random property whose probability can be quantified.

To formalize the above, we introduce some notation.
Given $\xi\in\Xi$ and $(g_1,\ldots,g_N)\in\calG^N$, we denote by $\Alg_\xi(g_1,\ldots,g_N)$ the solution of~\eqref{eq:optim-scenario} with $g_1,\ldots,g_N$.\footnote{Without loss of generality, we assume that the solution exists and is unique; see, e.g.,~\cite{calafiore2010random} for ways to handle problems with no solutions or non-unique solutions.
{\newgb Note that this definition of $A$ is very general: it includes any scenario optimization problems, and any method to obtain an exact \emph{or approximate} solution to it.
Allowing approximate solutions is particularly relevant when dealing with nonconvex or NP-hard problems.}}
Given $x\in\Xsf$ and $\xi\in\Xi$, we denote by $V_\Prob(x,\xi)$ its \emph{risk}, i.e., the probability of violating the constraint $g(x,\xi)\leq0$ with respect to $\Prob$: $V_\Prob(x,\xi)=\Prob[g(x,\xi)>0]$.
When $\Prob$ and $\xi$ are clear from the context, we write $\Alg$ and $V(x)$ for $\Alg_\xi$ and $V_\Prob(x,\xi)$.
The probability that the solution of~\eqref{eq:optim-scenario} is feasible for~\eqref{eq:optim-cc} is then given by
\begin{equation}\label{eq:PAC}
C_\xi(\epsilon,N)\coloneqq\Prob^N\!\left(\left\{ \gb\in\calG^N : V_\Prob(\Alg_\xi(\gb),\xi) \leq \epsilon \right\}\right),
\end{equation}
where $\gb$ is a shorthand notation for $(g_1,\ldots,g_N)$.
Several lower bounds on $C_\xi(\epsilon,N)$ have been proposed in the literature.
These bounds generally depend on an intrinsic quantity of the problem of interest, called its \emph{complexity}~\cite{campi2008theexact,alamo2009randomized,calafiore2010random,margellos2015ontheconnection}.
The definition of complexity varies from one approach to another.
For non-degenerate convex scenario programs, the dimension of the decision variable can be used as a complexity measure~\cite{campi2008theexact,calafiore2010random}.
Other common complexity measures include the VC dimension, the Rademacher complexity and the compression size~\cite{alamo2009randomized,margellos2015ontheconnection}.
Yet, sharp bounds on these quantities are generally elusive, resulting in overly conservative sample size requirements.


We address the problem of finding sharp upper bounds on $C_\xi(\epsilon,N)$ in the context of \emph{repetitive scenario design}, that is, when~\eqref{eq:optim-scenario} is solved repeatedly {\newgb with $\xi$ (or more precisely its effect on~\eqref{eq:optim-scenario}) varies slowly}.
The goal is that ultimately $N$ is close to the smallest value such that $C_\xi(\epsilon,N)\geq\beta$, where $\beta\in[0,1]$ is a given \emph{confidence parameter}.
In other words, the problem we address is the following:

\begin{problem}\label{prob:offline}
Let $\calG$, $\Prob$ and $\Alg$ be as above.
Given $\epsilon\in[0,1]$ and $\beta\in[0,1]$, let $N_\star(\xi)=\min\,\{N\in\N : C_\xi(\epsilon,N)\geq\beta\}$.
Find a repetitive scenario design algorithm such that, with high probability, $N_t\approx N_\star(\xi_t)$ for all $t\in\N_{>0}$ large enough, where $N_t$ is the sample size used at step $t$ of the algorithm {\newgb and $\xi_t$ is the value of the external parameter at step $t$}.
\end{problem}

We stress out that a small change in $N$ can have a large impact on the computational complexity of solving~\eqref{eq:optim-scenario}; for nonconvex problems with exact computations, the dependence in $N$ can be exponential or worse.

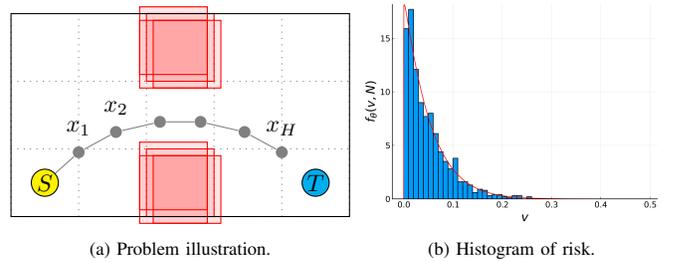
\begin{figure}[t]
    \centering
    \centering
    \begin{subfigure}[b]{0.53\linewidth}
        \centering
        \begin{tikzpicture}[scale=0.9,transform shape]
            \draw[step=1.0,gray,thin,dotted] (0,0) grid (5,3);
            \coordinate (S) at (0.5,0.5);
            \coordinate (T) at (4.5,0.5);
            \path[draw=red,fill=red,fill opacity=0.1] (2,0) rectangle (3,1);
            \path[draw=red,fill=red,fill opacity=0.1,shift={(-0.1,-0.1)}] (2,0) rectangle (3,1);
            \path[draw=red,fill=red,fill opacity=0.1,shift={(-0.1,+0.1)}] (2,0) rectangle (3,1);
            \path[draw=red,fill=red,fill opacity=0.1,shift={(+0.1,-0.1)}] (2,0) rectangle (3,1);
            \path[draw=red,fill=red,fill opacity=0.1] (2,3) rectangle (3,2);
            \path[draw=red,fill=red,fill opacity=0.1,shift={(-0.1,-0.1)}] (2,3) rectangle (3,2);
            \path[draw=red,fill=red,fill opacity=0.1,shift={(-0.1,+0.1)}] (2,3) rectangle (3,2);
            \path[draw=red,fill=red,fill opacity=0.1,shift={(+0.1,-0.1)}] (2,3) rectangle (3,2);
            \node[fill=gray,circle,minimum width=5pt,inner sep=0pt] (m1) at (1.0,0.95) {};
            \node[fill=gray,circle,minimum width=5pt,inner sep=0pt] (m2) at (1.55,1.25) {};
            \node[fill=gray,circle,minimum width=5pt,inner sep=0pt] (m3) at (2.2,1.4) {};
            \node[fill=gray,circle,minimum width=5pt,inner sep=0pt] (m4) at (2.8,1.4) {};
            \node[fill=gray,circle,minimum width=5pt,inner sep=0pt] (m5) at (3.45,1.25) {};
            \node[fill=gray,circle,minimum width=5pt,inner sep=0pt] (m6) at (4.0,0.95) {};
            \path[draw=gray] (S) -- (m1) -- (m2) -- (m3) -- (m4) -- (m5) -- (m6);
            \draw[black,fill=yellow] (S) circle (0.2) node[inner sep=0pt] (S) {$S$};
            \draw[black,fill=cyan] (T) circle (0.2) node[inner sep=0pt] (T) {$T$};
            \node[above=1pt of m1] {$x_1$};
            \node[above=1pt of m2] {$x_2$};
            \node[above=1pt of m6] {$x_H$};
            \draw[black] (0,0) rectangle (5,3);
        \end{tikzpicture}
        \caption{Problem illustration.}
        \label{fig:example-mpc-problem}
    \end{subfigure}%
    \hfill
    \begin{subfigure}[b]{0.45\linewidth}
        \centering
        \includegraphics[width=\linewidth]{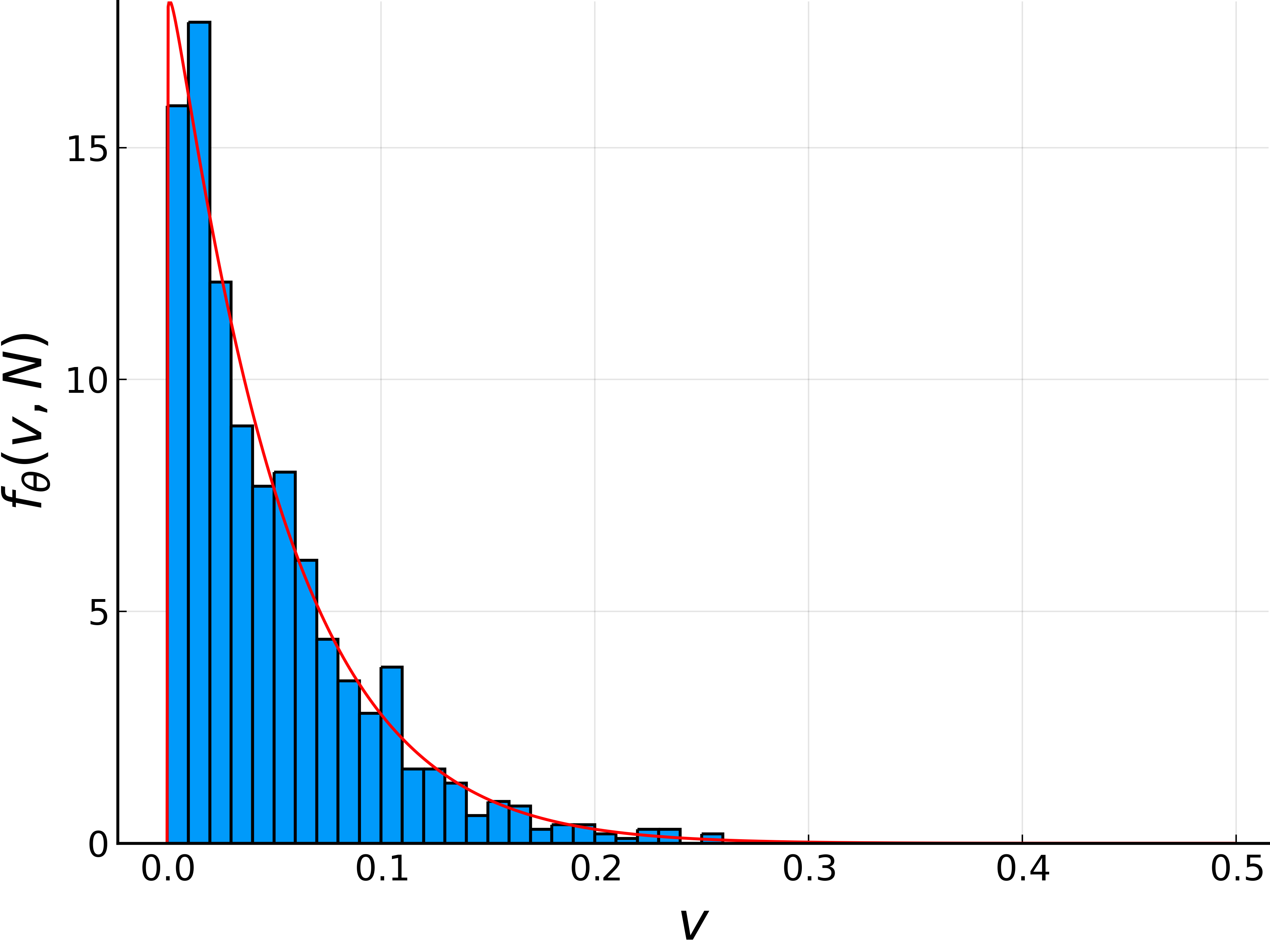}
        \caption{Histogram of risk.}
        \label{fig:example-mpc-risk}
    \end{subfigure}
    \caption{Path planning problem of Example~\ref{exa:path-planning}.
    {(a) The path is represented by a linked sequence of $H$ viapoints (gray dots).
    Each viapoint must avoid the obstacles (red regions) and two consecutive viapoints cannot be more than $\delta$ units of distance apart.
    (b) Histogram (using $M=1000$ experiments) of the risk of paths computed with $N=20$ sampled positions of the obstacles {\newgb(hence, a total of $20\,000$ sampled constraints)}, and $H=100$ viapoints.
    The red curve is $f_\theta(\cdot,N)$ where $\theta=\theta_\star(\calD)$ {\newgb and $\calD$ is the data set obtained from the $M$ experiments (cf.~Section~\ref{sec:approach})}.}}
    \label{fig:example-mpc}
\end{figure}

We conclude this section with an example:

\begin{example}\label{exa:path-planning}
Consider the problem of finding the shortest path between a source location $\xi_S$ and a target location $\xi_T$ while avoiding moving obstacles, as depicted in Fig.~\ref{fig:example-mpc-problem}.
In the spirit of~\cite{jankowski2023vpsto,brudermuller2024ccvpsto}, we parametrize the path by using $H$ ``viapoints'', denoted by $x_1,\ldots,x_H\in\R^2$, which are our decision variables.
Each viapoint needs to avoid the obstacles and two consecutive viapoints cannot be more than $\delta$ units of distance apart.
Problem~\ref{eq:optim-robust} becomes
\begin{equation}\label{eq:optim-mpc-robust}
\begin{array}{cl}
\min\limits_{x_1,\ldots,x_H\in\Bsf} & \lVert x_H-\xi_T\rVert_2 \\
\text{s.t.} & \lVert x_t - x_{t-1} \rVert_2 \leq \delta, \quad t\in[H], \\
& g(x_t,\xi) \leq 0, \quad t\in[H], \quad g\in\calG,
\end{array}
\end{equation}
where $x_0=\xi_S$, $\Bsf=[0,5]\times[0,3]$, and for each $g\in\calG$, $g(x,\xi)$ is the signed distance (positive in case of collision) between $x$ and the obstacles in one specific position.
In the example of Fig.~\ref{fig:example-mpc-problem}, each $g\in\calG$ has the form
\[
g(x,\xi; \xt_l,\xt_u) = \min\,\{\tfrac12-\lVert x-\xt_l\rVert_\infty,\tfrac12-\lVert x-\xt_u\rVert_\infty\},
\]
where $\xt_l,\xt_u\in\R^2$ are the centers of the lower and upper obstacles, respectively.
Note that in Fig.~\ref{fig:example-path} there is no path from $\xi_S$ to $\xi_T$ that satisfies the constraint $g(x)\leq0$ for all $\xt_l,\xt_u\in\R^2$.
Hence, the robust problem~\eqref{eq:optim-mpc-robust} is typically infeasible.
Assume that a probability distribution $\Prob$ on the value of $\xt_l$ and $\xt_u$, e.g., $[\xt_l;\xt_u]\sim\calN(\mu,\Sigma^2)$, is given.
In this case, the chance-constrained problem~\eqref{eq:optim-cc} becomes relevant and captures the fact that the path avoids the obstacles with high probability.
Solving the chance-constrained problem, which is often intractable, can be approached by solving the scenario problem~\eqref{eq:optim-scenario}.
In this case, it amounts to sample $N$ values of $(\xt_l,\xt_u)$ and solve the optimization problem
\begin{equation}\label{eq:optim-mpc-scenario}
\!\begin{array}{@{}cl@{}}
\min\limits_{x_1,\ldots,x_H\in\Bsf} & \lVert x_H-\xi_T\rVert_2 \\
\text{s.t.} & \lVert x_t - x_{t-1} \rVert_2 \leq \delta, \quad t\in[H], \\
& g(x_t,\xi;\xt_{l,i},\xt_{u,i}) \leq 0, \quad t\in[H], \; i\in[N].
\end{array}
\end{equation}
The determination of an optimal sample size $N$ is the central question of this paper.\hfill$\triangleleft$
\end{example}

{\newgb

In the rest of this paper, we assume that $\xi$ is fixed, and we address Problem~\ref{prob:offline} under this assumption.
The idea is that the results obtained for $\xi$ fixed will be useful in practice for problems with slowly-varying $\xi$ (as we showcase in numerical experiments in Section~\ref{sec:experiments}).
A theoretical analysis of the case of slowly-varying $\xi$ is left for future work.

}

\section{PROPOSED APPROACH}\label{sec:approach}

Our approach consists in learning a function $f_\theta(v,N)$ that approximates for each $N\in\N$ the pdf of $V(\Alg(g_1,\ldots,g_N))$ when $g_1,\ldots,g_N$ are sampled i.i.d.~from $\Prob$.
In other words, we aim to find a value of $\theta$ such that $f_\theta(v,N)\approx f(v,N)$, where $f(v,N) \coloneqq \frac{\dd}{\dd v} \Prob^N\!\left(\left\{ \gb\in\calG^N : V(\Alg(\gb)) \leq v \right\}\right)$.
To do that, given a data set $\calD=\{(v_j,N_j)\}_{j=1}^M\subseteq [0,1]\times\N$, we define the Maximum Likelihood Estimator (MLE):
\begin{equation}\label{eq:MLE}
\theta_\star(\calD)=\argmax_{\theta\in\Theta} \, \ell(\theta;\calD), 
\end{equation}
where $\ell(\theta;\calD) \coloneqq \frac1M \sum_{j=1}^M \log(f_\theta(v_j,N_j))$.


In the rest of this section, we first explain how we build the data set $\calD$ in the context of repetitive scenario design.
Then, we describe the parametrization that we consider for $f_\theta$ and explain how we compute $\theta_\star(\calD)$ for this parametrization.
Finally, we explain how we estimate the optimal sample size $N$ given a parameter $\theta$ satisfying $f_\theta(v,N)\approx f(v,N)$.

\subsection{Data collection}\label{ssec:data-collection}

In the context of repetitive scenario design,~\eqref{eq:optim-scenario} is solved repeatedly with different sample sets.
Let $N_t$ be the sample size used in step $t$, and let $x_t$ be the returned solution.
For each $t\in\N_{>0}$, let $v_t=V(x_t)$.
This gives at each step $t\in\N_{>0}$, the data set $\calD_t=\{(v_s,N_s)\}_{s=1}^t$, from which the MLE $\theta_t=\theta_\star(\calD_t)$ can be computed.

\begin{remark}
Given $x$, $V(x)$ can sometimes be easily computed; this is the case, e.g., for the problems in Sections~\ref{ssec:fixed-complexity} and~\ref{ssec:time-varying}.
Otherwise, it can always be easily approximated to high accuracy by using a Bernoulli test (e.g., the problems in Section~\ref{ssec:non-fixed-complexity}); see\ifextended\ Appendix~\ref{app:risk-evaluation}\else~\cite[Appendix~A]{berger2025online}\fi\ for details.
\end{remark}

\subsection{Distribution shape and fitting}

In this paper, we propose the following parametrization of $f_\theta$: for $v\in(0,1]$ and $N\in\N_{>\theta}$,
\begin{equation}\label{eq:beta-simple}
f_\theta(v,N) = \frac{v^{\theta-1}(1-v)^{N-\theta}}{B(\theta,N-\theta+1)}, \quad \theta\in\Theta\coloneqq\R_{>0},
\end{equation}
where $B$ is the Beta function (cf.~Section~\ref{sec:intro}).
Note that for all $N\in\N_{>\theta}$, $f_\theta(\cdot,N)$ is the pdf of the distribution $\mathrm{Beta}(\theta,N-\theta+1)$~\cite[\S25]{johnson1995continuous}.
The reasons for this choice of $f_\theta$ will be apparent in Theorem~\ref{thm:convergence}.
Note that \eqref{eq:beta-simple} requires that $v>0$ and $N>\theta$.
However, we want to allow that $\theta\geq N_j$ or $v_j=0$ for some $(v_j,N_j)\in\calD$.
Therefore, we extend the domain of $f_\theta$ as follows: for $v\in[0,1]$ and $N\in\N$,
\begin{equation}\label{eq:beta}
f_\theta(v,N) = \left\{
\begin{array}{@{\:}c@{\quad}l}
1 & \text{if} \; v=0 \; \text{or} \; N=0, \\
\frac{v^{\theta-1}(1-v)^{N-\theta}}{B(\theta,N-\theta+1)} & \text{if} \; v\neq0 \; \text{and} \; N>\theta, \\
N v^{N-1} & \text{otherwise.}
\end{array}\right.
\end{equation}

We discuss below some properties of $f_\theta$ that are useful for optimizing the log-likelihood $\ell(\cdot;\calD)$:

\begin{proposition}\label{prop:upper}
Let $v\in[0,1]$ and $N\in\N$.
The function $\theta\mapsto f_\theta(v,N)$ is upper semi-continuous on $\Theta$.
\end{proposition}

\ifextended%
\begin{proof}
If $v\neq1$, $f_\theta(v,N)$ is clearly continuous in $\theta$.
If $v=1$, $f_\theta(v,N)=0$ if $\theta<N$, and $f_\theta(v,N)=\max\{1,N\}$ if $\theta\geq N$, showing upper semi-continuity.
\end{proof}%
\fi

\begin{proposition}\label{prop:concave}
Let $v\in(0,1)$ and $N\in\N_{>0}$.
The function $\theta\mapsto\log(f_\theta(v,N))$ is smooth and strictly concave on $(0,N)$.
\end{proposition}

\ifextended%
\begin{proof}
Let $\theta\in(0,N)$.
It holds that
\begin{align*}
&\log(f_\theta(v,N)) = (\theta-1)\log(v)+(N-\theta)\log(1-v) \\
&\quad -\log(\Gamma(\theta))-\log(\Gamma(N-\theta+1))+\log(\Gamma(N+1)),
\end{align*}
where $\Gamma$ is the Gamma function (cf.~Section~\ref{sec:intro}), and we used that $B(a,b)=\frac{\Gamma(a)\Gamma(b)}{\Gamma(a+b)}$~\cite[\S6]{abramowitz1972handbook}.
Since $\log\circ\,\Gamma$ is smooth and strictly convex on $\R_{>0}$~\cite[\S6]{abramowitz1972handbook}, this concludes the proof.
\end{proof}%
\fi

\begin{corollary}\label{cor:smooth-mle}
Let $\calD=\{(v_j,N_j)\}_{j=1}^M\subseteq[0,1]\times\N$.
It holds that $\ell(\cdot;\calD)$ is upper semi-continuous on $\Theta$, and for every interval $I\subseteq\Theta\setminus\{N_j\}_{j=1}^M$, $\ell(\cdot;\calD)$ is either constant on $I$, or is smooth and strictly concave on $I$.
\end{corollary}

\ifextended%
\begin{proof}
First, fix $j\in[M]$.
By Proposition~\ref{prop:upper}, it holds that $\theta\mapsto\log(f_\theta(v_j,N_j))$ is upper semi-continuous on $\Theta$.
If $v_j=0$ or $N_j=0$, then $\theta\mapsto\log(f_\theta(v_j,N_j))$ is constant on $\Theta$.
If $v_j=1$ and $N_j>0$, then $\theta\mapsto\log(f_\theta(v_j,N_j))$ is constant ($-\infty$) on $(0,N_j)$ and is constant on $[N_j,\infty)$.
If $v_j\in(0,1)$ and $N_j>0$, then by Proposition~\ref{prop:concave}, $\theta\mapsto\log(f_\theta(v_j,N_j))$ is smooth and strictly concave on $(0,N_j)$ and is constant on $[N_j,\infty)$.
Finally, since $\ell(\theta;\calD) = \frac1M \sum_{j=1}^M \log(f_\theta(v_j,N_j))$, this concludes the proof.
\end{proof}%
\fi

Corollary~\ref{cor:smooth-mle} implies that the maximum of $\ell(\theta;\calD)$ exists and is finite.
It also implies that $\ell(\theta;\calD)$ can be maximized easily by considering separately each maximal interval $I$ of $\Theta\setminus\{N_j\}_{j=1}^M$; indeed, on each $I$, $\ell(\theta;\calD)$ is either constant or smooth and concave.
The maximum of a concave function defined on $\R$ can be computed very efficiently and reliably, e.g., using iterative methods.
In our numerical experiments, we used Newton--Raphson's method with an educated guess for the initial iterate; this method proved very efficient in all experiments, eliminating the need for additional measures or adjustments.

\begin{example}
The MLE $\theta_\star(\calD)$ for the problem in Example~\ref{exa:path-planning} is represented in Fig.~\ref{fig:example-mpc-risk}.
We collected $M=1000$ data points by solving the scenario problem $1000$ times with $N=20$ and $H=100$.\hfill$\triangleleft$
\end{example}

\subsection{Update of the sample size}

Given a parameter value $\theta$ such that $f_\theta(v,N)\approx f(v,N)$, a risk tolerance $\epsilon\in[0,1]$ and a confidence $\beta\in[0,1]$, we estimate the optimal sample size $N_\star$ in Problem~\ref{prob:offline} as follows.
By definition of $f$, it holds that $C(\epsilon,N)=\int_0^\epsilon f(v,N)\,\dd v$.
This gives the optimal sample size estimate
\[
\Nh(\theta;\epsilon,\beta) = \min \left\{ N\in\N : \int_0^\epsilon f_\theta(v,N) \, \dd v \geq \beta \right\}.
\]
Since the value of the integral is increasing with $N$, $\Nh(\theta;\epsilon,\beta)$ can be computed efficiently by using bisection.

The overall algorithm is presented in Algo.~\ref{algo:online-scenario}.
Note that if $V(x)$ is estimated by using a Bernoulli test as explained in\ifextended\ Appendix~\ref{app:risk-evaluation}\else~\cite[Appendix~A]{berger2025online}\fi, then $\Prob$ needs not to be known precisely by the algorithm; a \emph{generative model} (i.e., an oracle generating i.i.d.~samples from $\Prob$) is sufficient.

\begin{remark}
The parameter $\Nmax$ in Algo.~\ref{algo:online-scenario} can be arbitrarily large.
Its purpose is only to simplify the analysis of the convergence of the algorithm (see Theorem~\ref{thm:convergence}).
In future work, we will work on removing this parameter, even though in practice, $\{N_t\}_{t=1}^\infty$ will always be bounded due to hardware and software limitations.
\end{remark}

\newcommand\mycommfont[1]{\footnotesize\ttfamily#1}
\SetCommentSty{mycommfont}

\begin{algorithm}
\caption{Repetitive Scenario Design}
\label{algo:online-scenario}
\DontPrintSemicolon
\KwData{$\epsilon\in[0,1]$, $\beta\in[0,1]$, $N_1\in\N$, $\Nmax\in\N$.}
$\calD_0\gets\emptyset$ \tcp{Data set}
\For{$t=1,2,\ldots$}{
    Draw $N_t$ i.i.d.~samples $g_1,\ldots,g_{N_t}\sim\Prob$\;
    Let $x_t \gets \Alg(g_1,\ldots,g_{N_t})$\;
    Let $\calD_t \gets \calD_{t-1} \cup \{(V_\Prob(x_t),N_t)\}$\;
    Let $\theta_t \gets \theta_\star(\calD_t)$ \tcp{see~\eqref{eq:MLE}}
    Let $N_{t+1} \gets \min\,\{\Nh(\theta_t;\epsilon,\beta),\Nmax\}$\;
}
\end{algorithm}

\section{ANALYSIS OF THE ALGORITHM}\label{sec:theory}

We demonstrate the convergence of Algo.~\ref{algo:online-scenario} when applied to \emph{fixed-complexity} scenario problems, defined below.
This notion generalizes the notion of fully-supported convex scenario programs, which have been extensively studied in the scenario optimization literature~\cite{campi2008theexact,calafiore2010random,garatti2022risk,romao2021tight,romao2023ontheexact}.\footnote{\newgb It is worth noting that fixed-complexity scenario problems satisfy the non-degeneracy assumption~\cite{campi2008theexact,calafiore2010random,garatti2022risk}.}

\begin{definition}\label{def:fixed-complexity}
Given $\calG$, $\Prob$, $\Alg$ and $d\in\N_{>0}$, we say that $(\calG,\Prob,\Alg)$ has \emph{fixed complexity} $d$ if for all $N\in\N_{\geq d}$, the following holds with probability one on $(g_1,\ldots,g_N)\in\calG^N$: there is a unique subset $(g_{i_1},\ldots,g_{i_d})$ with $1\leq i_1<\ldots<i_d\leq N$ such that $\Alg(g_1,\ldots,g_N)=\Alg(g_{i_1},\ldots,g_{i_d})$.
\end{definition}

{\newgb

\begin{remark}
When applied to problems with non-fixed $\xi$, the analog condition would be that for all $\xi\in\Xi$, $(\calG,\Prob,\Alg_\xi)$ has fixed complexity $d_\xi$, where $d_\xi$ changes slowly over time.
A formal analysis is left for future work.
\end{remark}

}

When applied to fixed-complexity problems, Algo.~\ref{algo:online-scenario} converges to the optimal sample size, thereby providing a valid solution to Problem~\ref{prob:offline}.
This is formalized in the next theorem, which is the main theoretical result of this paper:

\begin{theorem}\label{thm:convergence}
Let $(\calG,\Prob,\Alg)$ have fixed complexity $d\in\N_{>0}$.
Let $0<\epsilon<\beta<1$ and $N_\star = \min \left\{ N\in\N : C(\epsilon,N) \geq \beta \right\}$.
Consider the sequences $\{N_t\}_{t=1}^\infty$ and $\{\theta_t\}_{t=1}^\infty$ generated by Algo.~\ref{algo:online-scenario}.
Assume that $N_\star\leq\Nmax$.
The following holds with probability one: (i) $\theta_t\to d$, and (ii) there is $t_0\in\N_{>0}$ such that for all $t\in\N_{\geq t_0}$, $N_t\in[N_\star-1,N_\star+1]$. 
\end{theorem}

\ifextended%
The proof relies on the following lemma, which is essentially the convergence of the MLE $\theta_\star$ for fixed $N$.

\begin{lemma}\label{lem:MLE}
Let $\theta_\circ>0$ and $N\in\N_{\geq\theta_\circ}$.
Let $\{v_t\}_{t=1}^\infty$ be a sequence of i.i.d.~random variables with distribution $\mathrm{Beta}(\theta_\circ,\allowbreak N-\theta_\circ+1)$.
For all $t\in\N_{>0}$, let $\calD_t=\{(N,v_j)\}_{j=1}^t$.
Define the function $\ell_N:\Theta\to\R$ by
\begin{align*}
&\ell_N(\theta) = (\thetat-1)\Psi(\theta_\circ)+(N-\thetat)\Psi(N-\theta_\circ+1) \\
&\quad -\log(\Gamma(\thetat))-\log(\Gamma(N-\thetat+1))+\log(\Gamma(N+1)) \\
&\quad -(N-1)\Psi(N+1),
\end{align*}
where $\thetat=\min\{\theta,N\}$, $\Gamma$ is the Gamma function and $\Psi$ is the Digamma function (cf.~Section~\ref{sec:intro}).
With probability one on $\{v_t\}_{t=1}^\infty\sim\Prob^\infty$, it holds that $\lim_{t\to\infty}\sup_{\theta\in\Theta}\,\lvert \ell(\theta;\calD_t)-\ell_N(\theta)\rvert=0$.
Also, $\ell_N$ is continuous, and if $N>\theta_\circ$, it has a unique maximizer at $\theta=\theta_\circ$.
\end{lemma}

The proofs of the lemma and the theorem are presented in Appendices~\ref{app:proof-lem-MLE} and~\ref{app:proof-thm-convergence}, respectively.%
\fi

The theoretical analysis is for the moment restricted to the class of fixed-complexity scenario problems.
However, in the next section, we show that our approach performs also very well in practice on a series of challenging repetitive scenario design problems, including non-fixed-complexity problems, nonconvex constraints and time-varying distributions.

\ifextended%
\begin{remark}
The study of \emph{consistency} of MLEs---i.e., convergence to the correct parameter value---is a classical topic in statistics and probability theory.
Classical results in the literature include convergence in probability under some mild assumptions, and almost sure convergence under additional assumptions such as concavity; see, e.g.,~\cite{newey1994large}.
These results were unfortunately not directly exploitable for our proof of Theorem~\ref{thm:convergence} because (i) the function $f_\theta$ is not concave in $\theta$ (it is even not continuous), and (ii) we needed a stronger result to prove the convergence of $\theta$ to $d$ when using samples from $\mathrm{Beta}(d,N-d+1)$ with \emph{different} values of $N$.
Therefore, we proved Lemma~\ref{lem:MLE}, which states a stronger result, namely the almost sure \emph{uniform} convergence of $\ell(\theta;\calD_t)$ to $\ell_N(\theta)$.
The proof relies on the particular shape of $f_\theta$.
\end{remark}%
\fi

\section{NUMERICAL EXPERIMENTS}\label{sec:experiments}

In all experiments, we used $\epsilon=0.1$ and $\beta=0.9$.\footnote{The codes are available at \url{https://github.com/guberger/OnlineScenarioOptimization.jl}.}

\subsection{Fixed-complexity problems}\label{ssec:fixed-complexity}

\begin{figure}
    \centering
    \includegraphics[width=\linewidth]{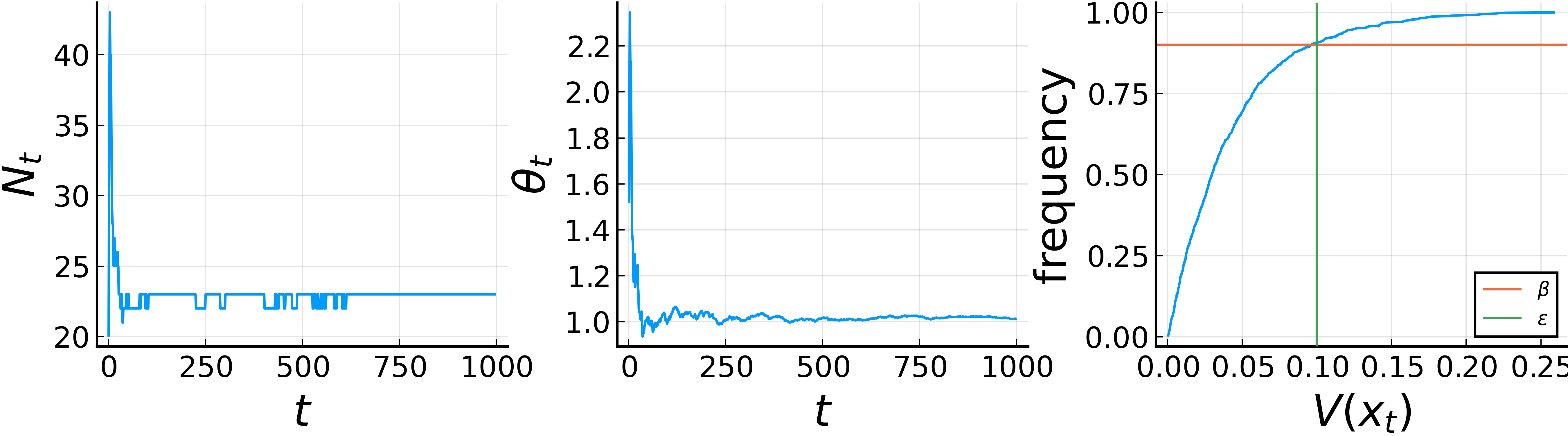}\\[2pt]
    \includegraphics[width=\linewidth]{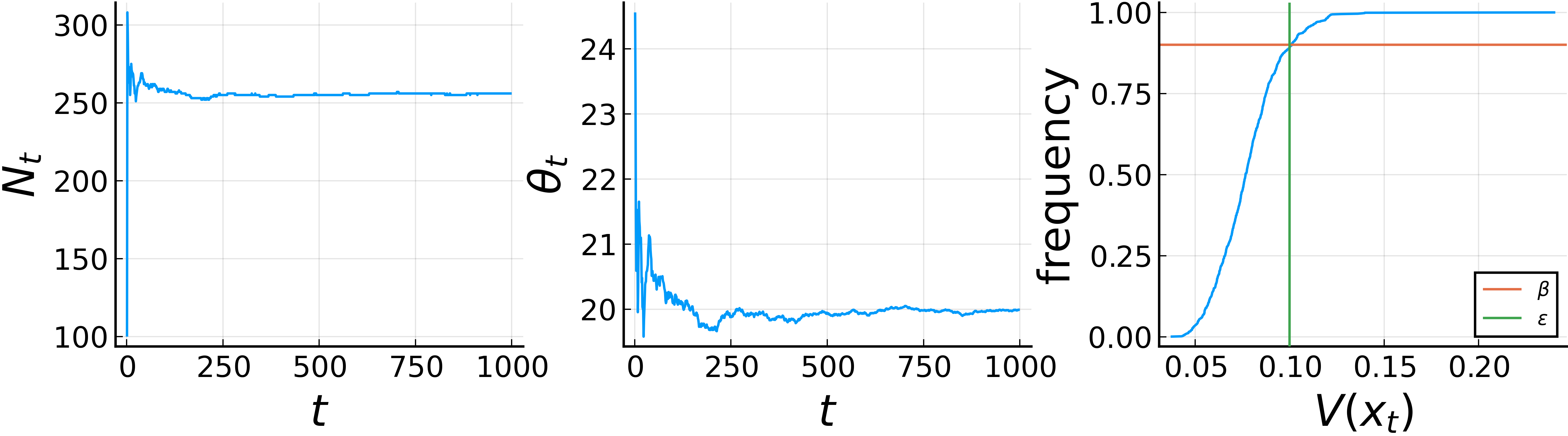}
    \caption{Fixed-complexity problems; top:~\eqref{eq:fixed-one}, bottom:~\eqref{eq:fixed-more}.
    Evolution of $N_t$ and $\theta_t$ over $T=1000$ steps and cumulative of $\{V(x_t)\}_{t=1}^T$.
    We observe that the algorithm converges very fast to $d$ and $N_\star$.
    We also observe that the violation probability is smaller than $0.1$, (almost) $90\%$ of the time.}
    \label{fig:fixed-complexity}
\end{figure}

We start with two fixed-complexity scenario problems.
The first problem consists in solving
\begin{equation}\label{eq:fixed-one}
\min_{x\in\R}\:x \quad \text{s.t.} \quad x\geq u, \quad u\sim\calN(1,2).
\end{equation}
This problem has fixed complexity $d=1$.
The evolution of $N_t$ and $\theta_t$ over the first $1000$ steps is represented in Fig.~\ref{fig:fixed-complexity}.
We also represented the cumulative of $V(x_t)$.
We notice that $\theta_t\to1$, and $V(x_t)\leq\epsilon$ with frequency at least $\beta$.

{\newgb

\begin{remark}
Note that for all our experiments, we reported the cumulative distribution of $\{V(x_t)\}_{t=1}^T$, where $T$ is the total number of steps.
This is because we are interested in the frequency with which $V(x_t)$ exceeds the risk tolerance over the \emph{whole} repetitive scenario design process.
\end{remark}

}

The second problem is
\begin{equation}\label{eq:fixed-more}
\min_{x\in\R^{20}}\: \sum_{i=1}^{20} x^{(i)} \quad \text{s.t.} \quad u^\top x\leq 1, \quad u\sim\calN(0,I).
\end{equation}
This problem has fixed complexity $d=20$.
The evolution of $N_t$ and $\theta_t$ over the first $1000$ steps is represented in Fig.~\ref{fig:fixed-complexity}.
We also represented the cumulative of $V(x_t)$.
We notice that $\theta_t\to20$, and $V(x_t)\leq\epsilon$ with frequency almost $\beta$.

\subsection{Beyond fixed-complexity problems}\label{ssec:non-fixed-complexity}

\begin{figure}
    \centering
    \includegraphics[width=\linewidth]{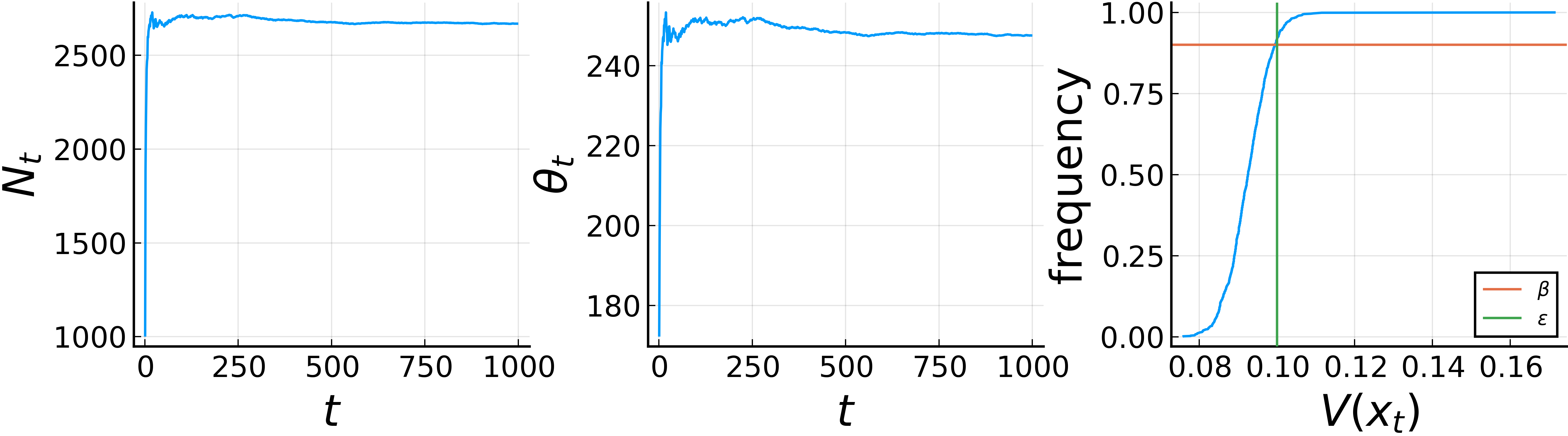}\\[2pt]
    \includegraphics[width=\linewidth]{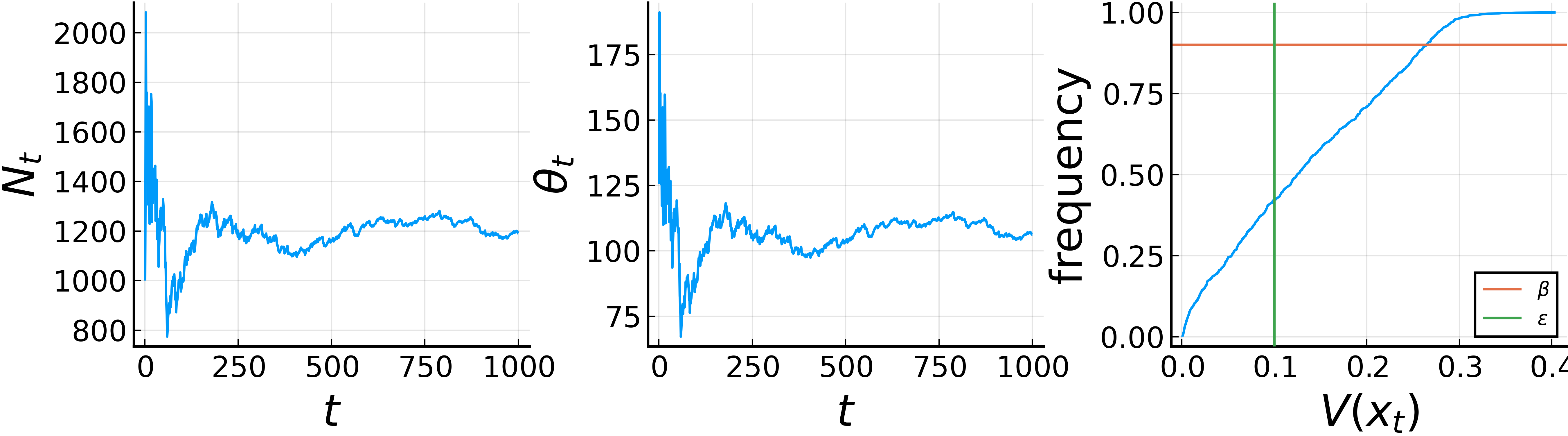}
    \caption{Non-fixed complexity problems from~\cite{garatti2022risk}; top:~\eqref{eq:garatti} with $\Prob$ as in~\cite[Fig.~8]{garatti2022risk}, bottom:~\eqref{eq:garatti} with $\Prob$ as in~\cite[Fig.~9]{garatti2022risk}.
    Evolution of $N_t$ and $\theta_t$ over $T=1000$ steps and cumulative of $\{V(x_t)\}_{t=1}^T$.
    On the top, we observe that the algorithm converges very fast to some $\theta_\circ$ and $N_\circ$, even though the complexity is not fixed.
    We also observe that the violation probability is smaller than $0.1$, $90\%$ of the time.
    On the bottom, however, the algorithm does not seem to converge to some $\theta$ or $N$, and the violation probability is smaller than $0.1$, less than $90\%$ of the time.
    An explanation for that is that the problem has ``degenerate complexity''.}
    \label{fig:garatti}
\end{figure}

Next, we consider more general scenario problems.
The first problem, borrowed from~\cite{garatti2022risk}, consists in solving
\begin{equation}\label{eq:garatti}
\min_{x\in\R^{400}}\: \sum_{i=1}^{400} x^{(i)} \quad \text{s.t.} \quad \min_{i=1}^{400}\, x^{(i)}\geq u, \quad u\sim\Prob.
\end{equation}
We let $\Prob$ be the probability distribution used in~\cite[Fig.~8]{garatti2022risk}.
The results are given in Fig.~\ref{fig:garatti}.
We notice that $V(x_t)\leq\epsilon$ with frequency at least $\beta$.

\ifextended\relax\else\addtolength{\textheight}{-1.8cm}\fi

We also tried with the probability distribution used in~\cite[Fig.~9]{garatti2022risk}.
The results are given in Fig.~\ref{fig:garatti}.
In this case, however, the frequency of $V(x_t)\leq\epsilon$ is \emph{smaller than} $\beta$.
This failure to meet the safety requirements can be explained by the fact that the problem is far from having fixed complexity (see~\cite[Fig.~9]{garatti2022risk}).
In future work, we plan to investigate algorithmic ways to detect such problems that have a non-fixed complexity with high variance, and provide sound ways to converge to their optimal sample size.

\begin{remark}
For these two problems, the risk was estimated using the technique in\ifextended\ Appendix~\ref{app:risk-evaluation}\else~\cite[Appendix~A]{berger2025online}\fi, with $S=10^4$, giving an accuracy of $\eta=0.025$ with probability $1-10^{-5}$.
\end{remark}

\begin{figure}
    \centering
    \includegraphics[width=\linewidth]{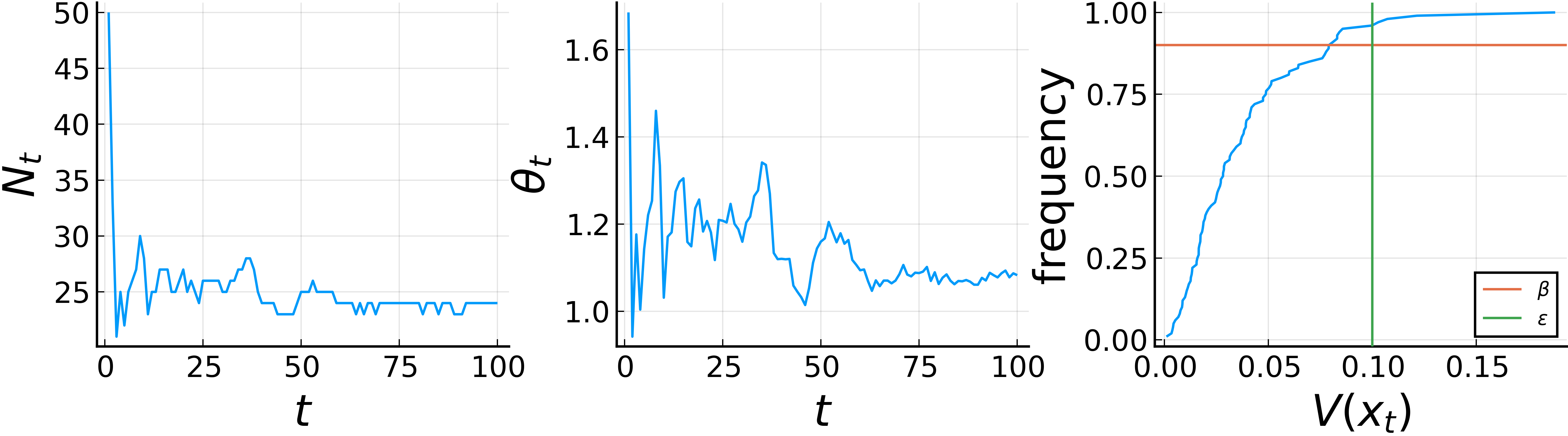}
    \includegraphics[width=\linewidth]{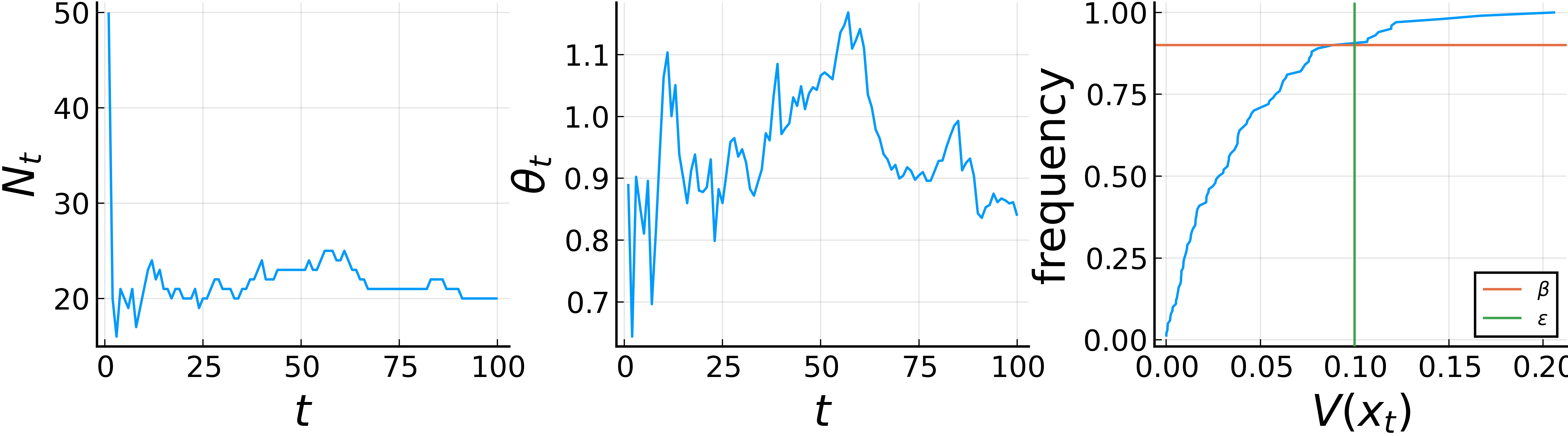}
    \caption{Path planning problem from Example~\ref{exa:path-planning}; top: steady distribution (see Section~\ref{ssec:non-fixed-complexity}), bottom: time-varying distribution (see Section~\ref{ssec:time-varying}).
    Evolution of $N_t$ and $\theta_t$ over $T=100$ steps and cumulative of $\{V(x_t)\}_{t=1}^T$.
    On the top, we observe that the algorithm converges fast to some $\theta_\circ$ and $N_\circ$, even though the complexity is not fixed.
    On the bottom, we observe that the value of $\theta_t$ and $N_t$ is time-varying, as is the distribution.
    For both, we observe that the violation probability is smaller than $0.1$, $90\%$ of the time.}
    \label{fig:path-planning}
\end{figure}

We also applied our technique on the path planning problem in Example~\ref{exa:path-planning} with $H=100$, $\delta=0.045$, $\xt_l=[\frac52,y-0.8]$ and $\xt_u=[\frac52,y+0.8]$, where $y\sim\calN(1.5, 0.05)$.
Note that the problem is nonconvex and the dimension of the decision variable is $200$.
The results are presented in Fig.~\ref{fig:path-planning}.
We notice that $V(x_t)\leq\epsilon$ with frequency at least $\beta$.

\subsection{Time-varying distribution}\label{ssec:time-varying}

Finally, we modify the path planning problem in Example~\ref{exa:path-planning} so that the distribution on the constraints is time-varying.
Namely, we let $\xt_l(t)=[\frac52+\sin(0.1t),y-0.3]$ and $\xt_u(t)=[\frac52+\sin(0.1t),y+0.3]$, where $y\sim\calN(1.5, 0.05)$ and $t\in\N_{>0}$ is the time step.
Because the distribution is shifting, we put more weight on the most recent data points with a rule proportional to the time step: $w_t=t$.
The results are presented in Fig.~\ref{fig:path-planning}.
We notice that although the distribution is shifting, $V(x_t)\leq\epsilon$ with frequency at least $\beta$.

\ifextended

\section{CONCLUSIONS}


We proposed an online-learning framework to learn the optimal sample size of scenario design problems, with application to repetitive scenario design.
We showed that our approach converges toward the optimal sample size for the class of fixed-complexity scenario design problems (which extend fully-supported convex scenario programs).
We also showed that it converges empirically toward the optimal sample size on a wide range of scenario design problems, including non-fixed-complexity problems, nonconvex constraints and time-varying distributions.

In future work, we plan to extend the formal convergence guarantees beyond fixed-complexity problems.
For that, we plan to use an adaptive re-weighting of the samples, putting more weight on the samples with larger risk.
The weighting will be updated until the probability of risk violation of the learned distribution corresponds to the one of the empirical distribution.
We also plan to remove the technical constraint $N_t\leq\Nmax$ in Algo.~\ref{algo:online-scenario}, which entails the (very mild) assumption that $N_\star\leq\Nmax$ in Theorem~\ref{thm:convergence}.

\section*{APPENDIX}

\renewcommand{\thesubsection}{\Alph{subsection}}

\subsection{Sample-based risk estimation and Bernoulli test}\label{app:risk-evaluation}

If computing the risk $V_\Prob(x)$ of a solution $x$ is too difficult computationally, one can resort to a statistical test to estimate $V_\Prob(x)$ with high accuracy.
The test consists in drawing $S$ i.i.d.~constraints $g_1,\ldots,g_S$ in $\calG$ according to $\Prob$ and computing the \emph{empirical risk}: $\vh(g_1,\ldots,g_S;x)=\frac1S \sum_{i=1}^S \oneb[g_i(x)>0]$.
For any error bound $\eta>0$, one can bound the probability that $V_\Prob(x)\notin[\vh(\gb;x)-\eta,\vh(\gb;x)+\eta]$ as
\[
\Prob^S\!\left(\left\{\gb\in\calG^S : \lvert \vh(\gb;x) - V_\Prob(x)\rvert > \eta \right\}\right) \leq 2e^{-2\eta^2S},
\]
by using Hoeffding's inequality~\cite[Lemma~B.6]{shalevshwartz2014understanding}.

\subsection{Proof of Lemma~\ref{lem:MLE}}\label{app:proof-lem-MLE}

The following lemma (see, e.g., \cite[\S25]{johnson1995continuous}) will be instrumental:

\begin{lemma}\label{lem:exp-beta}
Let $a,b>0$ and $X\sim\mathrm{Beta}(a,b)$.
It holds that $\Exp[\log(X)]=\Psi(a)-\Psi(a+b)$, where $\Psi$ is the Digamma function (cf.~Section~\ref{sec:intro}).
\end{lemma}

We proceed with the proof of Lemma~\ref{lem:MLE}.
Remember that $\{v_t\}_{t=1}^\infty$ is i.i.d.~with distribution $\mathrm{Beta}(\theta_\circ,N-\theta_\circ+1)$.
Since the Beta distribution is continuous, it holds with probability one that for all $t\in\N_{>0}$, $v_t\in(0,1)$.
Hence, with probability one, it holds that for all $t\in\N_{>0}$ and $\theta\in\Theta$,
\begin{align*}
&\ell(\theta;\calD_t) = \frac1t \sum_{s=1}^t (\thetat-1)\log(v_s)+(N-\thetat)\log(1-v_s) \\
&\quad -\log(\Gamma(\thetat))-\log(\Gamma(N-\thetat+1))+\log(\Gamma(N+1)).
\end{align*}
The Strong Law of Large Numbers~\cite[Theorem~8.2.7]{athreya2006measure} and Lemma~\ref{lem:exp-beta} imply that with probability one,
\begin{align*}
\textstyle \frac1t \sum_{s=1}^t \log(v_s) &\to \Psi(\theta_\circ)-\Psi(N+1), \\
\textstyle \frac1t \sum_{s=1}^t \log(1-v_s) &\to \Psi(N-\theta_\circ+1)-\Psi(N+1).
\end{align*}
This shows that $\sup_{\theta\in\Theta}\,\lvert\ell(\theta;\calD_t)-\ell_N(\theta)\rvert\to0$.
Finally, we show that $\ell_N$ has a unique maximizer at $\theta=\theta_\circ$.
On $(0,N]$, $\ell_N$ is smooth and strictly concave, with derivative
\[
\ell_N'(\theta)=\Psi(\theta_\circ)-\Psi(N-\theta_\circ+1)-\Psi(\theta)+\Psi(N-\theta+1).
\]
The derivative is zero at $\theta=\theta_\circ$.
Hence, with strict concavity, it follows that for all $\theta\in(0,N]\setminus\{\theta_\circ\}$, $\ell_N(\theta)<\ell_N(\theta_\circ)$.
On $[N,\infty)$, $\ell_N$ is constant and equal to $\ell_N(N)$.
This shows that if $N>\theta_\circ$, the unique maximizer of $\ell_N$ is $\theta=\theta_\circ$.\hfill\QED

\subsection{Proof of Theorem~\ref{thm:convergence}}\label{app:proof-thm-convergence}

The following result, adapted from~\cite[Theorem~3.3]{calafiore2010random}, will be instrumental in the proof:

\begin{lemma}\label{lem:fixed-beta}
Let $(\calG,\Prob,\Alg)$ have fixed complexity $d\in\N_{>0}$.
Let $N\in\N$.
(i) If $N\geq d$, then $V(\Alg(g_1,\ldots,g_N))$, with $g_1,\ldots,g_N\sim\Prob$, has distribution $\mathrm{Beta}(d,N-d+1)$.
(iii) If $N<d$, then $V(\Alg(g_1,\ldots,g_N))$, with $g_1,\ldots,g_N\sim\Prob$, is one with probability one.
\end{lemma}


We proceed with the proof of Theorem~\ref{thm:convergence}.

\subsubsection*{Outline of the proof}

Based on Lemma~\ref{lem:fixed-beta}, the idea of the proof is to use Lemma~\ref{lem:MLE} to show that with probability one, $\ell(\cdot,\calD_t)$ converges \emph{uniformly} toward a convex combination of the functions $\ell_N$ with $N\in\{d+1,\ldots,\Nmax\}$.
Since all these functions have a unique maximizer at $\theta=d$, we deduce that $\theta_t\coloneqq\theta_\star(\calD_t)$ converges toward $d$.
Since we need $N>d$ in the convex combination (to ensure unique maximizer), we divide the proof into two parts: first, we show (using Lemma~\ref{lem:MLE}) the result under the assumption that $N_t>d$ for all $t\in\N_{>0}$; then, we show (using the shape~\eqref{eq:beta} of $f_\theta$) that there is a finite number of $t\in\N_{>0}$ such that $N_t\leq d$.

\subsubsection*{Preliminaries}

The only source of randomness is the drawing of $g\sim\Prob$.
Hence, each run of Algo.~\ref{algo:online-scenario} is associated with a realization of the sequence $\gbt\coloneqq\{g_s\}_{s=1}^\infty$.
Therefore, $\{N_t\}_{t=1}^\infty$ and $\{\theta_t\}_{t=1}^\infty$ depend on $\gbt$, and will therefore be denoted by $\{N_t(\gbt)\}_{t=1}^\infty$ and $\{\theta_t(\gbt)\}_{t=1}^\infty$.
Our goal is to show that with probability one on $\gbt\sim\Prob^\infty$, $\theta_t(\gbt)\to d$.


\subsubsection*{Part 1}

First, we show it under the additional assumption that for all $\gbt\in\calG^\infty$ and $t\in\N_{>0}$, $N_t(\gbt)>d$.
We will show that with probability one on $\gbt\sim\Prob^\infty$, $\ell(\theta;\calD_t(\gbt))$ converges uniformly toward the convex hull of $\{\ell_N\}_{N=d+1}^{\Nmax}$, where $\ell_N$ is as in Lemma~\ref{lem:MLE}.
Therefore, for each $\gbt\in\calG^\infty$ and $N\in\{d+1,\ldots,\Nmax\}$, let $I_N=\lvert\{t\in\N_{>0} : N_t(\gbt)=N\}\rvert$, and for each $t\in\N_{>0}$, let $\alpha_{t,N}(\gbt)=\lvert I_N(\gbt) \cap [1,t]\rvert / t$ and $\calD_{t,N}(\gbt)=\{(N_s(\gbt),v_s(\gbt)) : s\in I_N(\gbt) \cap [1,t] \}$.
Observe that for all $\gbt\in\calG^\infty$ and $t\in\N_{>0}$, $\sum_{N=d+1}^{\Nmax}\! \alpha_{t,N}(\gbt)=1$ and
\[\textstyle
\ell(\theta;\calD_t(\gbt)) = \sum_{N=d+1}^{\Nmax} \alpha_{t,N}(\gbt) \ell(\theta;\calD_{t,N}(\gbt)).
\]
For each $\gbt\in\calG^\infty$ and $t\in\N_{>0}$, let us define the function $\ell_t(\cdot;\gbt)=\sum_{N=d+1}^{\Nmax} \alpha_{t,N}(\gbt) \ell_N$, which is in the convex hull of $\{\ell_N\}_{N=d+1}^{\Nmax}$.
Lemma~\ref{lem:MLE} and (i) in Lemma~\ref{lem:fixed-beta} imply that with probability one on $\gbt\sim\Prob^\infty$, $\lim_{t\to\infty}\sup_{\theta\in\Theta}\,\lvert\ell(\theta;\calD_t(\gbt))-\ell_t(\theta;\gbt)\rvert=0$.\footnote{Indeed, Lemma~\ref{lem:MLE} implies that with probability one on $\gbt\sim\Prob^\infty$, if $\lvert I_N(\gbt)\rvert=\infty$, then $\lim_{t\to\infty}\sup_{\theta\in\Theta}\,\lvert\ell(\theta;\calD_{t,N}(\gbt))-\ell_N(\theta)\rvert=0$.}
Finally, by Lemma~\ref{lem:MLE} again, for each $N\in\N_{>d}$, $\ell_N$ has a unique maximizer at $\theta=d$.
Since $\theta_t=\theta_\star(\calD_t(\gbt))$, this implies that with probability one on $\gbt\sim\Prob^\infty$, $\theta_t(\gbt)\to d$.

\subsubsection*{Part 2}

Next, we show that with probability one on $\gbt\sim\Prob^\infty$, $\{t\in\N_{>0} : N_t(\gbt) \leq d\}$ is finite.
For that, we first show that for each $N\in\N_{<d}$, with probability one on $\gbt\sim\Prob^\infty$, $\lvert\{t\in\N_{>0} : N_t(\gbt)=N\}\rvert\leq1$.
Fix $t\in\N_{>0}$, $u\in\N_{>t}$ and $N\in\N_{<d}$.
We will show that with probability one on $\gbt\sim\Prob^\infty$: if $N_t(\gbt)=N$, then $N_u(\gbt)>N$.
By (ii) in Lemma~\ref{lem:fixed-beta}, it holds with probability one on $g_1,\ldots,g_N\sim\Prob$ that if $N_t(\gbt)=N$, then $v_t(\gbt)=1$.
Hence, by definition of $f_\theta$, the following holds with probability one on $\gbt\sim\Prob^\infty$: if $N_t(\gbt)=N$, then $\theta_{u-1}(\gbt)\geq N$.
Indeed, fix $\gbt\in\calG^\infty$, and assume that $N_t(\gbt)=N$ and $v_t(\gbt)=1$.
Then, since $f_\theta(1,N')=0$ whenever $N'>\theta$, we have that $\ell(\theta;\calD_{u-1}(\gbt))=-\infty$ if $\theta<N$, implying that $\theta_{u-1}(\gbt)\geq N$.
Moreover, by definition of $\Nh$ and since $\beta>\epsilon$, it holds that $N_u(\gbt)>\min\{\theta_{u-1}(\gbt),\Nmax-1\}$.
Indeed, $\int_0^\epsilon f_\theta(v,N')\,\dd v=\epsilon^{\max\{1,N'\}}<\beta$ if $\theta\geq N'$.
Since $t$ and $u$ were arbitrary, this shows that with probability one on $\gbt\sim\Prob^\infty$, $\lvert\{t\in\N_{>0} : N_t(\gbt)=N\}\rvert\leq1$.

It remains to show that with probability one on $\gbt\sim\Prob^\infty$, $I_d(\gbt)\coloneqq\{t\in\N_{>0} : N_t(\gbt) = d\}$ is finite.
By Lemma~\ref{lem:MLE} and (i) in Lemma~\ref{lem:fixed-beta}, it holds with probability one on $\gbt\sim\Prob^\infty$ that if $\lvert I_d(\gbt)\rvert=\infty$, then $\theta_t(\gbt)\to[d,\infty)$.
Also, for all $\gbt\in\calG^\infty$, if $\theta_t(\gbt)\to[d,\infty)$, then $N_t(\gbt)\to[d+1,\Nmax]$ since $\epsilon<\beta$ (same argument as above).
This is a contradiction, showing that with probability one on $\gbt\sim\Prob^\infty$, $\{t\in\N_{>0} : N_t(\gbt) \leq d\}$ is finite.
This concludes the proof that with probability one on $\gbt\sim\Prob^\infty$, $\lim_{t\to\infty}\sup_{\theta\in\Theta}\,\lvert\ell(\theta;\calD_t(\gbt))-\ell_t(\theta;\gbt)\rvert=0$, where $\ell_t$ is defined in Part~1, in the general case.
Hence, with probability one on $\gbt\sim\Prob^\infty$, $\theta_t(\gbt)\to d$.

\subsubsection*{Conclusion}

Finally, observe that $\Nh$ is pseudo-continuous, meaning that for all $\theta\in\Theta$, there is $\eta>0$ such that for all $\theta'\in\Theta$, $\lvert\theta-\theta'\rvert\leq\eta$ implies $\lvert\Nh(\theta)-\Nh(\theta')\rvert\leq1$.
Hence, with probability one on $\gbt\sim\Prob^\infty$, $N_t(\gbt)\to[\Nh(d)-1,\Nh(d)+1]$.
Furthermore, by Lemma~\ref{lem:fixed-beta}, it holds that $\Nh(d)=N_\star$.
This concludes the proof of the theorem.\hfill\QED

\fi






\bibliographystyle{IEEEtran}
\bibliography{myrefs}

\end{document}